\documentclass{article}

\advance\textheight -\headheight
\advance\textheight -\headsep
\oddsidemargin\paperwidth
\advance\oddsidemargin -\textwidth
\divide\oddsidemargin2 
\ifdim\oddsidemargin<.5truein \oddsidemargin.5truein \fi
\advance\oddsidemargin -1truein
\evensidemargin\oddsidemargin
\topmargin\paperheight \advance\topmargin -\textheight
\advance\topmargin -\headheight \advance\topmargin -\headsep
\divide\topmargin2
\ifdim\topmargin<.5truein \topmargin.5truein \fi
\advance\topmargin -1.4truein\relax

\usepackage{amsmath,amsthm,amssymb,amscd,enumerate,mathtools,enumitem}
\usepackage{amsfonts}
\usepackage{rotating}
\usepackage{euscript}
\usepackage{pst-node}
\usepackage{epsfig,verbatim}
\usepackage{stackengine,scalerel}
\stackMath

\def\be{\begin{equation}}
\def\ee{\end{equation}}

\def\mult{{\rm mult}} 
\def\D{{\rm D}} 
 
\def\C{{\mathbb C}} 
\def\f{\EuScript}
 
\def\P{{\mathbb P}}

\def\T{{\mathbb T}}
\def\D{{\mathbb D}}
\def\R{{\mathbb R}}

\def\phi{{\varphi}}
\def\v{{\varepsilon}}

\def\supp{{\rm supp\,}}
\def\bp{\begin{proposition}}
\def\ep{\end{proposition}}

\def\bt{\begin{theorem}}
\def\et{\end{theorem}}
\def\br{\begin{remark}}
\def\er{\end{remark}}
\def\be{\begin{equation}}
\def\bee{\begin{equation*}}
\def\l{\label}

\def\ee{\end{equation}}
\def\eee{\end{equation*}}
\def\bl{\begin{lemma}}
\def\el{\end{lemma}}
\def\bc{\begin{corollary}}
\def\ec{\end{corollary}}
\def\pr{\noindent{\it Proof. }}

\def\bd{\begin{definition}}
\def\ed{\end{definition}}

\def\hat{\widehat}

\newtheorem{theorem}{Theorem}[section]
\newtheorem{lemma}[theorem]{Lemma}
\newtheorem{definition}[theorem]{Definition}
\newtheorem{corollary}[theorem]{Corollary}
\newtheorem{proposition}[theorem]{Proposition}
\newtheorem{problem}[theorem]{Problem}

\theoremstyle{definition}

\theoremstyle{definition}
\newtheorem{remark}[theorem]{Remark}
\def\bpr{\begin{problem}}
\def\epr{\end{problem}}

\mathtoolsset{showonlyrefs}

\providecommand{\keywords}[1]{\textbf{Keywords:} #1}

\begin{document}

\date{}

\title{On dessins d'enfants with equal supports}

\author{Fedor Pakovich \footnote{Department of Mathematics, Ben Gurion University of the Negev, Israel, e-mail: \texttt{pakovich@math.bgu.ac.il.}}}

\maketitle

\hspace*{\fill} \emph{To George Shabat, on the occasion of his 70th birthday}

\vspace{5mm}

\begin{abstract}
 For a Belyi function $\beta:\C\P^1\rightarrow \C\P^1$  
ramified only over the points $-1,1,\infty$, a corresponding ``dessin d'enfant'' $\f D_{\beta}$ is defined as 
the set $\beta^{-1}([-1,1])$ considered as a bi-colored  graph on the Riemann sphere whose 
white and black vertices are points of the sets $\beta^{-1}\{-1\}$ and $\beta^{-1}\{1\}$ correspondingly.  
Merely  the set $\beta^{-1}([-1,1])$ without a graph structure is called a support of $\f D_{\beta}$.  In this note, we solve the following problem: under what conditions different dessins $\f D_{\beta_1}$ and $\f D_{\beta_2}$ have equal supports? 
\end{abstract}


\section{Introduction}
One of  consequences of the ``dessins d'enfants'' theory (see e. g.  \cite{des2}, \cite{des1}, \cite{des3}) is the fact that any bi-colored connected graph $\Gamma$ on the Riemann sphere has a ``true form". This means that there exists a Belyi function $\beta:\C\P^1\rightarrow \C\P^1$  
ramified only over the points $-1,1,\infty$ such that the set $\f D_{\beta}=\beta^{-1}([-1,1])$,  
considered as a bi-colored graph whose 
white and black vertices are points of the sets $\beta^{-1}\{-1\}$ and $\beta^{-1}\{1\}$ correspondingly, represents $\Gamma$ in the following sense: there exists 
 an orientation preserving homeomorphism of the sphere $\phi$ such that $\phi(\Gamma)=\beta^{-1}([-1,1])$ and $\phi$ maps white and black vertices of $\Gamma$ to white and black vertices of  $\beta^{-1}([-1,1])$ correspondingly. The graph $\f D_{\beta}$ is called a {\it dessin d'enfant} or simply a {\it dessin} corresponding to the Belyi function $\beta.$ In this note, we always assume that considered Belyi functions are rational functions on the sphere ramified only over the points $-1,1,\infty$. Correspondingly, all   considered dessins are   spherical.

The geometry of plane trees (that is, of graphs without cycles) given in a true form has been studied in several publications (see  \cite{bi}, \cite{pom}, \cite{ko1}, \cite{ko2}). 
In particular, the remarkable result of Bishop (\cite{bi}) states that for any compact, connected set $K \subset \C$ and any $\v > 0$ there exists a Shabat polynomial (that is, a polynomial Belyi function) $P$   such that the set $ P^{-1}([-1, 1])$ approximates $K$ to within $\v$  in the Hausdorff metric. Thus, the sets $P^{-1}([-1,1])$ are dense in all planar continua.

In virtue of the fundamental correspondence between bi-colored graphs and  Belyi functions, the dessin $\f D_{\beta}$ 
defines the correspondent Belyi function $\beta$ up to equivalence. In this note, we address the following problem, which seems to be especially interesting 
in view of the theorem of Bishop: let $\f D_{\beta}$ be a dessin, up to what extent merely the set $\beta^{-1}([-1,1])$, without a graph structure, defines $\f D_{\beta}$? 
In other words, under what conditions on $\f D_{\beta_1}$ and $\f D_{\beta_1}$ the sets $\beta_1^{-1}([-1,1])$ and $\beta_2^{-1}([-1,1])$ coincide as subsets of $\C\P^1$? 
For a Belyi function $\beta$,  we call the set 
 $\beta^{-1}([-1,1])$ a {\it support} of $\f D_{\beta}$ and denote it by $\supp\{\f D_{\beta}\}$. 

To describe pairs of dessins $\f D_{\beta_1}$ and  $\f D_{\beta_2}$ such that  
\be \l{opo}  \supp\{\f D_{\beta_1}\}=\supp\{\f D_{\beta_2}\},\ee  
it is enough to restrict ourselves to the case where corresponding Belyi functions $\beta_1$ and $\beta_2$  satisfy  the condition 
$\C(\beta_1,\beta_2)=\C(z),$ that is, generate the whole field of rational functions $\C(z)$. 
Indeed, the L\"uroth theorem yields that for arbitrary rational functions $\beta_1$, $\beta_2$ there exist rational functions  $\widehat \beta_1$,  $\widehat \beta_2$ and $W$ such that 
\be \l{eq2} \beta_1=\widehat  \beta_{1}\circ  W, \ \ \ \  \beta_2=\widehat  \beta_{2}\circ   W\ee  and $\C(\widehat \beta_1,\widehat  \beta_2)=\C(z)$. Moreover, it easy to see that the equalities \eqref{eq2}  and 
\be \l{e1} \beta_1^{-1}([-1,1])=\beta_2^{-1}([-1,1])\ee 
 imply the equality \be \l{eq3} \widehat\beta_1^{-1}([-1,1])=\widehat\beta_2^{-1}([-1,1]),\ee while equalities \eqref{eq2} and \eqref{eq3} imply equality \eqref{e1}. Finally, if $\beta_1$, $\beta_2$ are Belyi functions, then the chain rule implies that $\widehat \beta_1$,  $\widehat \beta_2$  are also Belyi functions.

 Examples of different dessins with equal supports are provided by ``segments'' and ``circles". By definition, a {\it segment} is a dessin whose support is homeomorphic to a segment, and a {\it circle} is a dessin whose support is homeomorphic to a circle. Notice that segments and circles can be characterized as dessins  all vertices of which have valency one or two. 
It is easy to see that the Chebyshev polynomial $\pm T_n,$ $n\geq 1$, is a Belyi function
corresponding to a segment with $n$ edges,  while $\pm \frac{1}{2}\left(z^n+\frac{1}{z^n}\right)$, $n\geq 1$,   
is a Belyi function corresponding to  a circle with $2n$ edges. 
Since for any  $n\geq 1$ the  preimage 
$(\pm T_n)^{-1}([-1,1])$ is the segment $[-1,1]$, for any  Belyi functions $\beta_1$ and $\beta_2$ such that $\f D_{\beta_1}$ and $\f D_{\beta_2}$ are segments,  there exists a M\"obius transformation $\mu$ such that  \be \l{op} \supp\{\f D_{\beta_1}\}=\supp \{\f D_{\beta_2\circ \mu}\}.\ee Similarly, 
equality \eqref{op} holds for some  M\"obius transformation $\mu$, if 
$\f D_{\beta_1}$ and $\f D_{\beta_2}$ are circles, 
since for any $n\geq 1$ the preimage $\left(\pm \frac{1}{2}\left(z^n+\frac{1}{z^n}\right)\right)^{-1}([-1,1])$ is the unit circle. Further examples of dessins satisfying \eqref{opo} can be obtained by formulas \eqref{eq2}, where $\widehat \beta_1$ and   $\widehat \beta_2$ are Belyi functions corresponding to segments or circles, and  $W$ is a rational function with an appropriate branching ensuring that $\beta_1$ and $\beta_2$ are also Belyi functions. In particular, for any Belyi function $\beta$ the equalities  
$$\supp\{\f D_{\beta}\}=\supp\{\f D_{\pm T_n\circ \beta}\}, \ \ \ n\geq 1,$$ 
hold.  

In this note, we prove the following result.

\bt \l{t1} 
Let $\f D_{\beta_1}$ and $\f D_{\beta_2}$ be dessins  such that $\supp\{\f D_{\beta_1}\}=\supp\{\f D_{\beta_2}\}$ and $\C(\beta_1,\beta_2)=\C(z)$. Then either $\f D_{\beta_1}$ and $\f D_{\beta_2}$ are segments, or 
$\f D_{\beta_1}$ and $\f D_{\beta_2}$ are circles. 
\et

Our approach to the problem is based on the observation that dessins d'enfants on the Riemann sphere 
are subsets of rational lemniscates. We recall that a {\it rational    lemniscate} is a 
curve in $\C$ defined by the equality 
\be \l{le} \f L_{\beta}=\{z\in \C: \, \vert \beta(z)\vert =1\},\ee where 
$\beta$ is a non-constant complex rational function. 
The geometry of lemniscates  is a classic subject of study. For example, Hilbert proved (\cite{hil}) that 
for every topological annulus $A \subset \C,$ there exists a polynomial $P$ whose
lemniscate is an analytic Jordan curve separating the two boundary components of $A.$
Lemniscates of entire and meromorphic functions were studied by Valiron (\cite{v}) and Cartwright (\cite{c}). 
For more recent results, we refer the reader to the papers \cite{bil}, \cite{havi}, \cite{ere},  
\cite{o},  \cite{op}, \cite{ps}, \cite{mum}, \cite{ste}, \cite{y}  and the bibliography therein. 
 
Since the real axis can be transformed to the unite circle by a M\"obius transformation $\mu$, for every rational function $\beta$ the set $\beta^{-1}([-1,1])$ is a subset of the lemniscate $\f L_{\mu\circ \beta}$.  Therefore, results about the geometry of lemniscates can be used for studying the geometry of dessins d'enfants. In particular, our proof of Theo\-rem \ref{t1} relies on the results of the recent papers  \cite{op}, \cite{ps}.

\section{Lemniscates and dessins d'enfants}
\subsection{Intersections of lemniscates}
Let us recall that a {\it finite  Blaschke product}  is a rational function of the form 
$$B=c\prod_{i=1}^m\frac{z-a_i}{1-\overline{a_i}z},$$
where  $\vert c \vert =1$ and $a_i$, $1\leq i \leq m,$ belong to the open unit disk $\D$.  
Correspondingly, a {\it quotient of finite Blaschke products} is a rational function of the form $B=B_1/B_2$, where $B_1$ and $B_2$ are  finite  Blaschke products. 
Notice that quotients of finite Blaschke products 
 can be characterized 
as rational functions that  map   
 the unit circle $\T$ 
to itself.

Let $P_1$ and $P_2$ be non-constant complex rational functions of degrees $n_1$ and $n_2$.  
It was shown in the recent paper \cite{ps} that 
the system of equations \be \l{sys1} \vert P_1(z)\vert =\vert P_2(z)\vert =1\ee has at most $(n_1+n_2)^2$ 
 solutions in $z\in \C$, unless 
\be \l{sys2}
 P_1=B_1\circ U, \qquad P_2=B_2\circ U
\ee
 for some rational functions $B_1$, $B_2$, and $U$, where  $B_1$ and $B_2$ are 
  quotients of finite Blaschke products. 
Geometrically, solutions of system \eqref{sys1} can be viewed as intersection points of the lemniscates $\f L_{P_1}$ and $\f L_{P_2},$ and 
in the forthcoming paper \cite{op} a sharp bound on the number of such points was obtained. Specifically, the result of \cite{op} states that, unless the condition \eqref{sys2} holds, the lemniscates $\f L_{P_1}$ and 
$\f L_{P_2}$ 
have at most $2n_1n_2$ intersection points, and this bound is the best possible.  
 
The results of \cite{op} and  \cite{ps}
imply the following statement.

\bt \l{t2} 
Let  $\beta_1,\beta_2$ be non-constant  rational functions, and $K_1$, $K_2$  subsets of $\R$ such that the intersection $\beta_1^{-1}(K_1)\cap\beta_2^{-1}(K_2)$ is infinite. Then 
there exist a rational function $W$ and rational functions with real coefficients $P_1$ and $P_2$ such that 
the equalities  
$ \beta_1=P_{1}\circ W$ and $ \beta_2=P_{2}\circ W$ hold. 
\et 
\pr 
Let us denote by $\hat\R$ the projectively extended real line,  $\hat\R=\R\cup \infty$, and  by $C$ the Cayley transformation, $$C(z)=\frac{z+i}{z-i}.$$
Since $C$ maps $\hat\R$  to  $\T$, if the intersection $\beta_1^{-1}(K_1)\cap\beta_2^{-1}(K_2)$ is infinite, then the  intersection $\f L_{C\circ \beta_1}\cap \f L_{C\circ \beta_2}$ 
is also infinite. Therefore, by the results of \cite{op} and \cite{ps},  there exist rational functions $B_1$, $B_2$, and $U$, where  $B_1$ and $B_2$ are 
  quotients of finite Blaschke products, such that 
\be \l{cl} 
 C\circ \beta_1=B_1\circ U, \ \ \   C\circ \beta_2=B_2\circ U.  
\ee 

Clearly, equalities \eqref{cl} imply the equalities 
\be \l{sys3}
\beta_1=(C^{-1}\circ  B_1 \circ C) \circ (C^{-1}\circ U), \ \ \ 
\beta_2=(C^{-1}\circ  B_2 \circ C) \circ (C^{-1}\circ U). 
\ee
On the other hand, since $C$ maps $\hat\R$  to  $\T$, a rational function   
$B$ is a quotient of a finite Blaschke products if and only if the rational function 
$P=C^{-1}\circ B\circ C$ maps $\hat\R$  to $\hat\R$. 
In turn, the last condition is equivalent 
 to the condition that $P$ has real coefficients (since  $P(z)-\overline{P}(z)=0$ for infinitely many $z\in \R$). Thus, \eqref{sys3} implies that the conclusion of the theorem holds for  $$ P_1=C^{-1}\circ B_1 \circ C,\ \ \ \  P_2=C^{-1}\circ B_2 \circ C,\ \ \ \ W=C^{-1}\circ U.\eqno{\Box}$$

\subsection{Proof of Theorem \ref{t1}}

 We deduce 
Theorem \ref{t1} from Theorem \ref{t2} and the following statement. 

\bt \l{p1} Let $P_1$ and $P_2$ be non-constant  rational functions with real coefficients such that $\R(P_1,P_2)=\R(z)$. 
Then $P_1^{-1}(\hat\R)\cap P_2^{-1}(\hat\R)=\hat\R\cup A$, where $A$ is a finite subset of $\C\P^1.$
\et
\pr 
Since $P_1$ and $P_2$ have real coefficients, $P_1(\hat\R)$ and $P_2(\hat\R)$ are subsets of $\hat\R,$ 
implying that 
$$\hat\R\subseteq P_1^{-1}(\hat\R)\cap P_2^{-1}(\hat\R).$$
On the other hand, since $\R(P_1,P_2)=\R(z)$, there exist non-zero polynomials with real coefficients $R_1(x,y)$ and $R_2(x,y)$   
such that \be \l{th} z=\frac{R_1(P_1(z),P_2(z))}{R_2(P_1(z),P_2(z))}.\ee Since for every point $z_0\in P_1^{-1}(\R)\cap P_2^{-1}(\R)$ the values $P_1(z_0)$ and $P_2(z_0)$ belong to $\R$, equality \eqref{th} implies that if $z_0$ belongs $P_1^{-1}(\R)\cap P_2^{-1}(\R)$, then 
$z_0$ belongs to $\R$, unless $z_0$ is a zero of $R_2(P_1(z),P_2(z))$. 
 This yields that the set 
$$(P_1^{-1}(\hat\R)\cap P_2^{-1}(\hat\R))\setminus\hat\R$$ is finite.  \qed 

\vskip 0.2cm 
\noindent{\it Proof of Theorem \ref{t1}.} 
Since the set  $\beta_1^{-1}([-1,1])$ is  infinite, it follows from \eqref{e1}  by 
 Theorem \ref{t2} taking into account the condition  $\C(\beta_1,\beta_2)=\C(z)$ 
that 
there exist a M\"obius transformation $\mu$ 
 and Belyi functions with real coefficients $\widetilde\beta_1$ and $\widetilde\beta_2$ such that 
 \be \l{e2} \beta_1=\widetilde\beta_{1}\circ \mu, \ \ \ \beta_2=\widetilde\beta_{2}\circ \mu,\ee  
 and  $\R(\widetilde\beta_1,\widetilde\beta_2)=\R(z)$.

For a rational function $\beta$ and a point $z\in \C\P^1$, we denote by $\mult_z \beta$ the multiplicity of $\beta$ at $z$. Let us observe that  for every $z_1\in \widetilde\beta_1^{-1}(\{-1,1\})$ and $z_2\in \widetilde\beta_2^{-1}(\{-1,1\})$ the inequalities 
$\mult_{z_1}\widetilde\beta_1\leq 2$ and $\mult_{z_2}\widetilde\beta_2\leq 2$ hold. 
Indeed, in a neighborhood of any point $z\in \widetilde\beta_1^{-1}(\{-1,1\})$  with $\mult_z \widetilde\beta_1=k\geq 3$  the set 
 \be \l{e3} T=\widetilde\beta_1^{-1}([-1,1])=\widetilde\beta_2^{-1}([-1,1])\ee 
 is homeomorphic to a  ``star'' with $k\geq 3$ branches. On the other hand, since \be \l{e4} T\subseteq \widetilde\beta_1^{-1}(\hat\R)\cap \widetilde\beta_2^{-1}(\hat\R),\ee  Theorem  \ref{p1} implies that 
$T\subset \hat\R\cup A$, where $A$ is a finite subset of $\C\P^1.$ The contradiction obtained shows that $\mult_{z_1} \widetilde\beta_1\leq 2$  for any $z_1\in \widetilde\beta_1^{-1}(\{-1,1\}).$ Similarly,  $\mult_{z_2} \widetilde\beta_2\leq 2$ for any $z_2\in \widetilde\beta_2^{-1}(\{-1,1\}).$ It follows from these inequalities  that 
 $\f D_{\beta_1}$ and $\f D_{\beta_2}$ are either segments or circles. Moreover, \eqref{e1}  implies that either the 
both dessins $\f D_{\beta_1}$ and $\f D_{\beta_2}$ are segments, or the both  are circles. 
 \qed

\begin{remark} 
Notice that the L\"uroth theorem implies that the condition \be \l{sat} \C(\beta_1,\beta_2)=\C(z)\ee is always satisfied if the degrees of $\beta_1$ and $\beta_2$ are coprime. Since $$T_{n_1n_2}=T_{n_2}\circ T_{n_1}, \ \ \ n_1,n_2\geq 1,$$ this implies in particular 
that the functions 
\be \l{beta0} \beta_1=\pm T_{n_1}, \ \ \ \beta_2=\pm T_{n_2}\ee satisfy \eqref{sat}  
 if and only if $\gcd(n_1,n_2)=1.$  On the other hand, the formulas
$$\pm  \frac{1}{2}\left(z^n+\frac{1}{z^n}\right)=\pm T_n\circ \frac{1}{2}\left(z+\frac{1}{z}\right), \ \ n\geq 1,$$ show   that  the functions 
\be \l{beta} \beta_1=\pm  \frac{1}{2}\left(z^{n_1}+\frac{1}{z^{n_1}}\right), \ \ \ \beta_2=\pm  \frac{1}{2}\left(z^{n_2}+\frac{1}{z^{n_2}}\right),\ \ \ n_1,n_2\geq 1,\ee never satisfy \eqref{sat}, since both these functions are rational functions in $\frac{1}{2}\left(z+\frac{1}{z}\right).$ 

Nevertheless, 
changing $\beta_1$ and $\beta_2$ in \eqref{beta} to the functions 
$\beta_1$ and $\beta_2 \circ \mu,$ where $\mu$ is a convenient M\"obius transformation that maps $\T$ to itself,  one can obtain  functions satisfying \eqref{opo} and  \eqref{sat} for arbitrary $n_1$ and $n_2$ (a similar remark is applied to the functions defined by \eqref{beta0}). 
Indeed, if equalities \eqref{eq2} hold for some rational function  $W$ of degree at least two, then the chain rule implies  that the functions  $\beta_1$ and $\beta_2$ have common critical points. Thus, it is enough to prove that for any rational functions 
$ \beta_1$ and $\beta_2 $ there exists $\mu$ as above such that $\beta_1$ and $\beta_2 \circ \mu$ have no common critical points. To prove it, we remark first that since  M\"obius transformations that map $\T$ to itself have the form $\frac{az+b}{\bar b z +\bar a},$ where $a,b$ are complex numbers satisfying $\vert a\vert \neq \vert b\vert,$ one can find  $\mu$ such that zero and infinity are not critical points of $\beta_2 \circ \mu$. As soon as this condition is met, obviously one can find a rotation $z\rightarrow cz$, $\vert c\vert =1$, such that the functions 
$\beta_1$ and $\beta_2 \circ \mu\circ cz$ have no common critical points.

\end{remark}

\end{document}